# NEW ASSOCIATED CURVES $k$ – PRINCIPLE DIRECTION CURVES AND $N_k$ – SLANT HELIX

Çağla RAMİS [1*], Beyhan UZUNOĞLU [2*] AND Yusuf YAYLI [3*]
cramis@ankara.edu.tr [1], buzunoglu@ankara.edu.tr [2], yusuf.yayli@science.ankara.edu.tr [3]
* Department of Mathematics, Faculty of Science, Ankara University, Ankara, TURKEY

**ABSTRACT.** The Frenet frame is generally known an orthonormal vector frame for curves. But, it does not always meet the needs of curve characterizations. In this study, with the help of associated curves of any spatial curve we obtained a new orthonormal frame which has the property that the second vector makes a constant angle with a fixed direction for the spatial curve. Then, the curve is named as the $N_k$ – slant helix and the special conditions are obtained effectively.

# 1. Introduction

In the practical world, curves arise in many different disciplines, for instance as the decorative in art, the path of particles in physics, the profile of technical object in medicine, engineering etc. Mathematically, a curve can be defined easily as a continuos mapping from an interval $I \subset \mathrm{R}$ to $\mathrm{R}^n$. If $n=3,$ then the calling turns as the spatial curve. Unfortunately, the continuity is so weak for definition of properties of curve. In addition to continuity, if any spatial parametrized curve $\gamma(s)$ requires that $\gamma$ is $n$ times continuously differentiable, then $\gamma$ is called $n$ – differentiable curve or $\gamma \in C^n$. And $\gamma(s)$ is a regular parametrized curve on the interval $I$, when $\gamma' \neq 0$. In three-dimensional Euclidean space, to explain the geometric structure of any regular spatial curve, an orthonormal basis $\{T, N, B\}$ called the Frenet frame at each point of the curve is described by the formulas which are independently discovered by Jean Frédéric Frenet, in his thesis of 1847, and Joseph Alfred Serret in 1851 [4,11]. Since the tangent vector $T$ is known and $B = T \times N,$ the moving frame $\{T, N, B\}$ is uniquely determined by the principle normal vector $N$. Thus $N$ is known as the reference vector of moving frame. Moreover, different space curves are only distinguished by the way in which they bend and twist and they are quantitatively measured by the invariants called the curvature $\kappa$ and the torsion $\tau$ of a curve, respectively. The fundamental theorem of curves asserts that there exists a unique regular parametrized curve $\gamma$ which has specific curvature and torsion with a solution of differential equation, also D. J. Struik give the following form of curve with Taylor expansion

$$\gamma(s) = \gamma(0) + \left(s - \frac{s^3 \kappa^2}{6}\right)T + \left(\frac{s^2 \kappa}{2} + \frac{s^3 \kappa'}{6}\right)N + \frac{s^3 \kappa \tau}{6}B + o(s^3) \qquad (1.1)$$

Where $T, N, B$ are the Frenet vector fields of $\gamma$ [7]. Also, with the help of curvatures some special curve characterizations are given, for example a classical result stated by M. A. Lancret in 1802 and first proved by B. de Saint Venant in 1845 for helices and slant helix representation is obtained by S. Izumiya and N. Takeuchi in 2004 [6,8,12].

On the other hand, the Frenet frame does not sufficiently satisfy for all curve to obtain characterizations and therefore made some changes on it to find featured different orthonormal frame on curves, such as the famous of them are the Bishop frame and RMF [2,5,14]. Each frame has its own advantages, of course, in mathematics and the other disciplines like that the useful material RMF for computer graphics and animations, motion design, robotics e.g.

In this study, with the help of a new orthonormal frame we give a kind of helix characterization for a spatial curve $\gamma \in C^{n+2}$ which has a useless Frenet frame in terms of generalization. The new frame is obtained under a process which is implemented the integral curves of the reference vector of $\gamma$ called $k-$principle direction curve of $\gamma$ until any principle direction curve is a slant helix. Thus, we explained the type of any spatial curve $\gamma$ as a $N_k-$slant helix in terms of the orthonormal frame of the slant helix $\gamma_k$ the $k-$principle direction curve of $\gamma$. The axis of $\gamma$ is got in terms of the Darboux vector $W_k$ and the normal vector $N_k$, so the curve $\gamma$ is named as a $W_k-$Darboux helix, too. Moreover, we obtained the special subfamily of $N_k-$slant helices called curves of $N_k-$constant precession which have constant speed Darboux vector. Finally, a useful example is given for the condition $k=1,2$ and in conclusion section we mentioned that $N_k-$slant helix includes some curves in [1,10,13].

## 2. Preliminaries

Let $\gamma$ be a unit-parametrized $2-$differentiable regular curve that $\|\gamma'\|=1$ and two step derivatives of $\gamma, \gamma'$ and $\gamma''$ exist in $\mathbb{E}^3$. The associated Frenet frame of $\gamma$, on each point of the curve, is an orthonormal basis at that point. The orthonormal basis is composed of three vectors called tangent, principal normal and binormal unit vectors at that point and defined as

$$T = \gamma', \ N = \frac{\gamma''}{\|\gamma''\|}, \ B = T \times N \qquad (2.1)$$

respectively.

At each point of the curve $\gamma$ the planes spanned by $\{T,N\}$, $\{T,B\}$, and $\{N,B\}$ are called the osculating plane, the rectifying plane and the normal plane, respectively. Also derivatives of the Frenet frame satisfy the following equations

$$\begin{aligned} T' &= \kappa N \\ N' &= -\kappa T + \tau B \\ B' &= -\tau N \end{aligned} \qquad (2.2)$$

Where $\kappa$ is the curvature and $\tau$ is the torsion of the curve $\gamma$, in sense of differential theory they describe completely how the frame evolves in time along the curve. Thus, the geometry of a unit-speed curve depends only on the values of curvature and torsion:

**(a)** if $\kappa = 0$, we obtain a uniformly-parametrized straight line;

**(b)** if $\tau = 0$, $\gamma$ be a planar curve;

**(c)** if $\gamma$ is the unit speed curve with constant curvature and torsion zero, then $\gamma$ is part of a circle of radius $1/\kappa$.

**(d)** $\gamma$ is a helix which has the property that the tangent line at any point makes a constant angle with a fixed line called the axis if and only if the rate $\frac{\kappa}{\tau}$ is constant;

**(e)** $\gamma$ is a slant helix whose principle normal vector at any point on $\gamma$ makes a constant angle with a fixed line if and only if the geodesic curvature of $N$, $\sigma = \dfrac{\kappa^2 \left(\frac{\tau}{\kappa}\right)'}{\left(\kappa^2 + \tau^2\right)^{3/2}}$ is a constant function, etc.

Let a point $p$ move along the regular curve $\gamma(s)$, then the Frenet frame $T, N, B$ changes, so the spherical indicatrix occur. As the point moves along the curve, the Frenet frame makes an instant helical motion at each $s$ moment along an axis called Darboux axis of $\gamma(s)$. And the vector which indicates the direction of axis is called Darboux vector of $\gamma$ and expressed as

$$W = \tau T + \kappa B \tag{2.3}$$

and it has the following symmetrical properties:

$$\begin{aligned} W \times T &= T' \\ W \times N &= N' \\ W \times B &= B' \end{aligned} \tag{2.4}$$

Moreover, the Darboux vector provides a concise way of interpreting curvature $\kappa$ and torsion $\tau$ geometrically: curvature is the measure of the rotation of the Frenet frame about the binormal unit vector $B$, whereas torsion is the measure of the rotation of the Frenet frame about the tangent vector $T$.

Also one of the application area of special vector fields of $\gamma$ is the equation of the axis helices. If $\gamma$ is a helix or slant helix, then the axis is

$$u = \cos\theta T + \sin\theta B \quad \text{or} \quad u = \sin\theta W + \cos\theta N \tag{2.5}$$

where $\theta$ is the constant angle with the tangent or the normal vector field of $\gamma$, respectively. (see [9])

## 3. $k$ – Principal Direction Curve and $N_k$ – Slant Helix

In this section, we investigate the integral curves of $(n+2)$ – differentiable space curve $\gamma$ in terms of $\{T, N, B\}$. Also generalization of one of them called the principle direction curves and slant helices is given. Moreover substantial relations are obtained between characterizations of curves with the help of frames.

Now, we define some associated curves of a curve $\gamma$ in $\mathbb{E}^3$ with using integral of its spherical indicatrix $T(s), N(s)$ and $B(s)$.

**Definition 1.** *Let $\gamma(s)$ be a regular unit-speed curve in terms of $\{T, N, B\}$. The integral curves of $T(s), N(s)$ and $B(s)$ are called the tangent direction curve, principal direction curve and binormal direction curve of $\gamma$, respectively [3]* ·

According to the sense of Frenet frame, the tangent direction curve $\gamma_0 = \int T(s)ds = \gamma + \vec{c}$ which is the translation form of $\gamma$ along the vector $\vec{c}$ has the same frame with $\gamma$ as $\{T_0 = T,\ N_0 = N,\ B_0 = B\}$. Also the binormal direction curve has the permutation frame $\{B, N, T\}$. On the other hand the principal direction curve of $\gamma$, $\gamma_1 = \int N(s)ds$ has a new frame as,

$$T_1 = N$$

$$N_1 = \frac{N'}{\|N'\|} = \frac{-\kappa T + \tau B}{\sqrt{\kappa^2 + \tau^2}} \qquad (3.1)$$

$$B_1 = T_1 \times N_1 = \frac{\tau T + \kappa B}{\sqrt{\kappa^2 + \tau^2}}$$

where the tangent vector and the binormal vector of $\gamma_1$ are the principle normal vector and the unit Darboux vector of $\gamma$, respectively. Also curvatures of $\gamma_1$ are

$$\kappa_1 = \sqrt{\kappa^2 + \tau^2}, \quad \tau_1 = \sigma \kappa_1 \qquad (3.2)$$

where $\sigma$ is the geodesic curvature of $N_1$.

If the principal direction curve of $\gamma$ is a regular curve, then its principal-direction curve can be defined and the process repeated again for itself. So if $\gamma$ is a $(n+2)$ – differentiable curve, then the principle direction curves generated such as:

The first principle direction curve,

$$\gamma_1 = \int N(s)ds \text{ with the frame } \left\{ T_1 = N, N_1 = \frac{N'}{\|N'\|}, B_1 = T_1 \times N_1 \right\}$$

The second principle direction curve

$$\gamma_2 = \int N_1(s)ds \text{ with the frame } \left\{ T_2 = N_1, N_2 = \frac{N_1'}{\|N_1'\|}, B_2 = T_2 \times N_2 \right\}$$

$$\vdots$$

The $(n)$ th principle direction curve

$$\gamma_n = \int N_{n-1}(s)ds \text{ with the frame } \left\{ T_n = N_{n-1}, N_n = \frac{N_{n-1}'}{\|N_{n-1}'\|}, B_n = T_n \times N_n \right\}.$$

**Definition 2.** *Let $\gamma(s)$ be a unit-parametrized and $(n+2)$ – differentiable curve in terms of $\{T, N, B\}$ in $\mathbb{E}^3$ and $\gamma_0$ be the tangent direction curve of $\gamma$. The $k$ – principal direction curve of $\gamma$ is defined as*

$$\gamma_k(s) = \int N_{k-1}ds, \quad 1 \leq k \leq n \qquad (3.3)$$

*where $N_{k-1}$ is the principle normal vector of $\gamma_{k-1}$. Also $\gamma$ is called the main curve of $\gamma_k$.*

The Frenet frame and curvatures of $\gamma_k$ are

$$T_k = N_{k-1}, \; N_k = \frac{N_{k-1}'}{\|N_{k-1}'\|}, \; B_k = T_k \times N_k \qquad (3.4)$$

$$\kappa_k = \sqrt{\kappa_{k-1}^2 + \tau_{k-1}^2}, \quad \tau_k = \sigma_{k-1} \kappa_k \qquad (3.5)$$

where $N_{k-1}, \kappa_{k-1}, \tau_{k-1}, \sigma_{k-1}$ are the principle normal vector, curvature, torsion and geodesic curvature of $N_{k-1}$, respectively and $B_k$ is the unit Darboux vector of $\gamma_{k-1}$.
The Frenet equations are satisfied

$$\begin{pmatrix} T_k \\ N_k \\ B_k \end{pmatrix}' = \begin{pmatrix} 0 & \kappa_k & 0 \\ -\kappa_k & 0 & \tau_k \\ 0 & -\tau_k & 0 \end{pmatrix} \begin{pmatrix} T_k \\ N_k \\ B_k \end{pmatrix} \qquad (3.6)$$

and the Darboux vector of $\gamma_k$ is

$$W_k = \tau_k T_k + \kappa_k B_k \qquad (3.7)$$

$$\begin{aligned} W_k \times T_k &= T_k' \\ W_k \times N_k &= N_k' \\ W_k \times B_k &= B_k' \end{aligned} \qquad (3.8)$$

also the unit Darboux vector of $\gamma_k$ is denoted as $\overline{W}_k$, $1 \le k \le n$.
Moreover with the help of helix theory, we can talk about relation between $k$-principle direction curve and its main curve. And the following definition are given to keep the generality.

**Definition 3.** Let $\gamma$ be a $(n+2)$-differentiable spatial curve, $\gamma_0$ and $\gamma_k$ be the tangent direction curve and $k$-principal direction curves of $\gamma$, $1 \le k \le n$, respectively. The curve $\gamma$ is called $N_k$-slant helix which has the property that the principal normal vector of $\gamma_k$ makes a constant angle with a fixed line called the axis.

**Remark 1.** If the normal vector of any curve makes a constant angle with a fixed line, then we called the curve as a slant helix in the preliminaries section of the paper. In this theory, being slant helix is the case of $k = 0$.

After giving all materials which are the above-mentioned definitions, now let give the following theorems and results to complete the theory.

**Theorem 1.** Let $\gamma$ be a $(n+2)$-differentiable curve in $\mathbb{E}^3$, and $\gamma_k$ be the $k$-principal direction curves of $\gamma$, $1 \le k \le n$. The curve $\gamma$ is a $N_k$-slant helix if and only if the $k$-principal direction curve $\gamma_k$ is a slant helix.

**Proof.** From definition 4, if $\gamma$ is a $N_k$-slant helix, then he principal normal vector of $\gamma_k$ makes a constant angle with a fixed direction $\vec{u}$

$$\langle N_k, \vec{u} \rangle = \cos \theta.$$

So, this definition coincides with being the slant helix of $\gamma_k$.
This completes the proof.

**Corollary 1.** The main curve $\gamma$ is a $N_k$-slant helix if and only if $\sigma_k = \dfrac{\kappa_k^2 \left( \dfrac{\tau_k}{\kappa_k} \right)'}{\left( \kappa_k^2 + \tau_k^2 \right)^{3/2}}$ which is the geodesic curvature of $N_k$ is a constant function.

**Theorem 2.** Let $\gamma$ be a $(n+2)$-differentiable spatial curve, $\gamma_k$ be the $k$-principle direction curve of $\gamma$, $1 \le k \le n$. The curve $\gamma$ is a $N_k$-slant helix if and only if the tangent indicatrix of $\gamma_k$ is a helix.

**Proof.** Let $\gamma$ be a $N_k$-slant helix, then from the theorem 1 the $k$-principle direction curve $\gamma_k$

is a slant helix. So, its tangent indicatrix $T_k$ is a helix.

Conversely, if $T_k$ is a helix, then the rate its curvature $\kappa_T$ and torsion $\tau_T$ is constant. This means that,

$$\frac{\tau_T}{\kappa_T} = \frac{\kappa_k^2 \left(\frac{\tau_k}{\kappa_k}\right)'}{\left(\kappa_k^2 + \tau_k^2\right)^{3/2}} = \text{constant}$$

where $\tau_k, \kappa_k$ are the torsion and curvature of $\gamma_k$, respectively.
This completes the proof.

The Frenet frame of $\gamma_k$, $\{T_k, N_k, B_k\}$ is an orthonormal basis of $\mathbb{E}^3$ and this frame is parametrized by $s$ which is the parameter of $\gamma$, so it can be located on the base curve of $\gamma$ and the other integral curves too. With the help of the frame sense, let define the axis of $N_k$ – slant helix.

**Theorem 3.** *Let $\gamma$ the main curve of $\gamma_k$ be a $N_k$ – slant helix in $\mathbb{E}^3$. The axis of $\gamma$ is*

$$\vec{u} = \sin\theta \overline{W}_k + \cos\theta N_k \tag{3.9}$$

*where $N_k, \overline{W}_k$ are the principal normal vector, the unit Darboux vector of $\gamma_k$ and $\theta$ is the constant angle with $N_k$.*

**Proof.** Let the main curve $\gamma$ be a $N_k$ – slant helix with the axis $\vec{u}$, then

$$\langle N_k, \vec{u} \rangle = \cos\theta \tag{3.10}$$

where $N_k$ is the principle normal vector of $\gamma_k$.
After the twice derivatives of Eq. (3.10),

$$\begin{aligned} -\kappa_k \langle T_k, \vec{u} \rangle + \tau_k \langle B_k, \vec{u} \rangle &= 0, \\ -\kappa_k' \langle T_k, \vec{u} \rangle + \tau_k' \langle B_k, \vec{u} \rangle &= \left(\kappa_k^2 + \tau_k^2\right)\cos\theta \end{aligned} \tag{3.11}$$

then the solutions of Eq. (3.11) are

$$\begin{aligned} \langle T_k, \vec{u} \rangle &= \tau_k \sigma_k \cos\theta, \\ \langle B_k, \vec{u} \rangle &= \kappa_k \sigma_k \cos\theta. \end{aligned} \tag{3.12}$$

So the axis of $\gamma$ is obtained as

$$\begin{aligned} \vec{u} &= \left(\tau_k \sigma_k T_k + N_k + \kappa_k \sigma_k B_k\right)\cos\theta \\ &= \sin\theta \overline{W}_k + \cos\theta N_k. \end{aligned} \tag{3.13}$$

**Corollary 2.** *If any spatial curv $\gamma$ is a $N_k$ – slant helix with the axis $\vec{u} = \sin\theta \overline{W}_k + \cos\theta N_k$, then it is obviously seen that $\overline{W}_k$ which is the unit Darboux vector of the $k$ – principle direction curve of $\gamma$ makes a constant angle with a fixed direction $\vec{u}$. Thus, the curve $\gamma$ is also called $\overline{W}_k$ – Darboux helix in 3-dimensional Euclidean space.*

# 4. Curves of $N_k$ – Constant Precession

In this section, we define a special subfamily of $N_k$ – slant helix called $N_k$ – constant precession in $\mathbb{E}^3$. Characteristic features of $N_k$ – constant procession are researched in terms of $k$ – principle direction curves.

Any spatial curve called curve of constant precession has the property that its Darboux vector revolves about a fixed direction $\vec{u}$ called the axis with constant angle and constant speed. As a consequence, its Frenet frame precesses about the axis, while its principal normal revolves about the axis with constant complementary angle and constant speed [10]. So, the family of constant precession is a kind of slant helix which has constant-length Darboux vector. With the help of constant precession theory, the following definition is given.

**Definition 4.** *Let $\gamma$ be a $(n+2)$ – differentiable spatial curve, $\gamma_0$ and $\gamma_k$ be the tangent direction curve and $k$ – principal direction curves of $\gamma$, $1 \leq k \leq n$, respectively. The curve $\gamma$ is called curve of $N_k$ – constant procession which has the property that the Darboux vector of $\gamma_k$ makes a constant angle with a fixed direction $\vec{u}$ and constant speed in terms of $\{T_k, N_k, B_k\}$.*

**Corollary 3.** *If $\gamma$ is a $N_k$ – slant helix with the Darboux vector $\parallel W_k \parallel = \omega$ (constant) in terms of $\{T_k, N_k, B_k\}$ then $\gamma$ is a curve of $N_k$ – constant precession.*

In previous section, we showed that if the main curve $\gamma$ is a $N_k$ – slant helix, then its axis lies in the plane spanned by $\overline{W}_k$, $N_k$ which are the unit Darboux vector and the principle normal vector of $\gamma_k$ in Eq. (3.9). If $\gamma$ is a curve of $N_k$ – constant procession with $\kappa_k^2 + \tau_k^2 = \omega^2$, then the axis of $\gamma$ is obtained such as

$$\vec{u} = W_k + \mu N_k. \tag{4.1}$$

Now, lets specialize the Darboux vector for the curve of $N_k$ – constant procession in the following theorem.

**Theorem 4.** *Any spatial curve $\gamma$ is a curve of $N_k$ – constant procession with the axis $\vec{u} = W_k + \mu N_k$, a necessary and sufficient condition is that*

$$\kappa_k = \omega \sin \mu \theta \quad \text{and} \quad \tau_k = \omega \cos \mu \theta \tag{4.2}$$

*where $\omega$ is the length of Darboux vector $W_k$.*

**Proof.** Let $\gamma$ be a curve of $N_k$ – constant procession with the axis $\vec{u} = W_k + \mu N_k$ with $\|W_k\| = \omega$, $\|\vec{u}\| = (\omega^2 + \mu^2)^{1/2}$ where $\omega$ and $\mu$ are constant. The differentiate of $\vec{u}$

$$\vec{u}' = (\tau_k' - \mu \kappa_k) T_k + (\kappa_k' + \mu \tau_k) B \tag{4.3}$$

is zero, therefore

$$\tau_k' - \mu \kappa_k = 0 \quad \text{and} \quad \kappa_k' + \mu \tau_k = 0 \tag{4.4}$$

The solution of this differential equation is uniquely

$$\kappa_k = \omega \sin \mu \theta \quad \text{and} \quad \tau_k = \omega \cos \mu \theta \tag{4.5}$$

This completes the proof.

**Example 1.** *Let $\gamma(s)$ be a unit speed curve with curvatures $\kappa(s) = \sin s \cos(\sin s)$, $\tau(s) = \sin s \sin(\sin s)$ and the curve $\gamma$ specifies neither helix nor slant helix in terms of the Frenet frame $\{T(s), N(s), B(s)\}$. Now lets obtain principle direction curves of $\gamma$ and characterize itself. Using the equations (3.5), the curvature and torsion of $\gamma_1 = \int N(s) ds$ the first principle direction*

curve of $\gamma$ are obtained as
$$\kappa_1 = \sin s \text{ and } \tau_1 = \cos s$$
So from Theorem 1 in [10], $\gamma_1$ is a curve of constant precession with the equation $\gamma_1(s) = (a_1(s), a_2(s), a_3(s))$

$$a_1(s) = \frac{(\sqrt{2}+1)^2}{2\sqrt{2}} \sin(\sqrt{2}-1)s - \frac{(\sqrt{2}-1)^2}{2\sqrt{2}} \sin(\sqrt{2}+1)s,$$

$$a_2(s) = -\frac{(\sqrt{2}+1)^2}{2\sqrt{2}} \cos(\sqrt{2}-1)s + \frac{(\sqrt{2}-1)^2}{2\sqrt{2}} \cos(\sqrt{2}+1)s,$$

$$a_3(s) = \frac{1}{\sqrt{2}} \sin s,$$

and the tangent vector field $T_1 = (t_1, t_2, t_3)$ is

$$t_1 = \frac{\sqrt{2}+1}{2\sqrt{2}} \cos(\sqrt{2}-1)s - \frac{\sqrt{2}-1}{2\sqrt{2}} \cos(\sqrt{2}+1)s,$$

$$t_2 = \frac{\sqrt{2}+1}{2\sqrt{2}} \sin(\sqrt{2}-1)s - \frac{\sqrt{2}-1}{2\sqrt{2}} \sin(\sqrt{2}+1)s,$$

$$t_3 = \frac{1}{\sqrt{2}} \cos s,$$

the principal vector field $N_1 = (n_1, n_2, n_3)$ is

$$n_1 = \frac{1}{\sqrt{2}} \cos \sqrt{2}s,$$

$$n_2 = \frac{1}{\sqrt{2}} \sin \sqrt{2}s,$$

$$n_3 = \frac{-1}{\sqrt{2}},$$

the binormal vector field $B_1 = (b_1, b_2, b_3)\}$ is

$$b_1 = -\sin\sqrt{2}s \cos s + \frac{1}{\sqrt{2}} \cos\sqrt{2}s \sin s,$$

$$b_2 = \cos\sqrt{2}s \cos s + \frac{1}{\sqrt{2}} \sin\sqrt{2}s \sin s,$$

$$b_3 = \frac{1}{\sqrt{2}} \sin s,$$

where $T_1 = N(s)$, $B_1$ is unit Darboux vector field of $\gamma$ also the Darboux vector of $\gamma_1$ is

$$W_1 = (\frac{1}{\sqrt{2}} \cos\sqrt{2}s, \frac{1}{\sqrt{2}} \sin\sqrt{2}s, \frac{1}{\sqrt{2}}).$$

From Theorem 1, the first principal direction curve of $\gamma$ is a slant helix so $\gamma$ is a $N_1$ – slant helix. Moreover the principal curve $\gamma_1$ is a curve of constant precession with constant speed Darboux, then $\gamma$ is called curve of $N_1$ – constant procession.

Because, $T_1(s)$ tangent indicatrix of $\gamma_1$ is a helix with the curvatures are $\bar{\kappa} = \frac{1}{\sin s}$, $\bar{\tau} = \frac{-1}{\sin s}$ and the fixed angle between tangent and the fixed direction is $\theta = \frac{-\pi}{4}$. So the axis of $\gamma$ is

$$\vec{u} = \sin\theta \overline{W}_1 + \cos\theta N_1$$
$$= (0, 0, -1).$$

*The second principal direction curve* $\gamma_2 = \int N_1 ds = \left(\frac{1}{2}\sin\sqrt{2}s, \frac{-1}{2}\cos\sqrt{2}s, \frac{-1}{\sqrt{2}}s\right)$ *is helix with curvatures and axis* $\kappa_2 = 1,\ \tau_2 = -1,\ \vec{v} = (0,0,1),$ *respectively. On the other hand the set of helices are the special subfamily of slant helices. So the curve* $\gamma$ *is a* $N_2$ – *slant helix with the axis* $\vec{v}$ *in terms of* $\{T_2, N_2, B_2\}$ *the Frenet frame of* $\gamma_2$.

*With using the Eq. (1.1) and the mathematica programme, we obtain approximately figures of the curve γ, γ₁ and γ₂.*

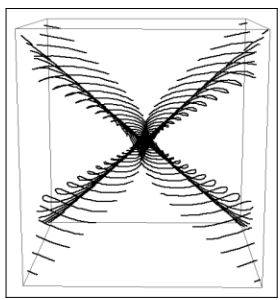

*Figure 1. The main curve γ(s)*

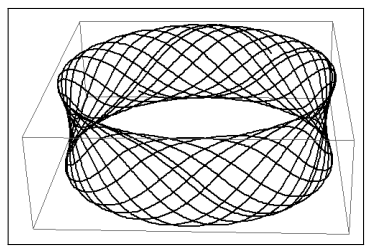

*Figure 2. The first principle direction curve γ₁*

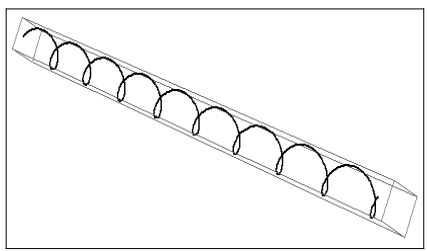

*Figure 3. The second principle direction curve γ₂*

## 5. Conclusion

Obtaining identification of curves is one of important study area in mathematics. Generally, the Frenet frame and curvatures of curves are used, the most famous characterizations are given for helix and slant helix by Lancret, Izumiya and Takeuchi [6,8]. Also Scofield calls the slant helices with constant speed Darboux vector field as curves of constant precession [10]. However the Frenet frame is not enough all curve, then new orthonormal frames are built on the curve to make up the disadvantages of it. In this study, new orthonormal frames which can be positioned on the regular spatial curve $\gamma$ are obtained to give characterization for $\gamma$ with the help of principle direction curves of $\gamma$. So, we called the curve $\gamma$ as a $N_k-$slant helix, if $\gamma_k$ the $k-$principle direction curve of $\gamma$ is a slant helix in terms of the frame $\{T_k, N_k, B_k\}$. And showed that the axis of $N_k-$slant helix lies in the plane $Sp\{\overline{W}_k, N_k\}$ where $\overline{W}_k$ and $N_k$ are the unit Darboux and the principle normal vector of $\gamma_k$, thus $\gamma$ is also named as a $\overline{W}_k-$Darboux helix in three dimensional Euclidean space. When the condition is $k=1$, the curve $\gamma$ is a $N_1-$slant helix or $C-$slant helix according to Uzunoğlu, Gök and Yaylı [13]. Moreover, with the help of the spherical image of tangent indicatrix and its own parametrized Frenet frame, a new special curve called $k-$slant helix is introduced by Ali [1], but $k-$slant helix and $N_k-$slant helix totally specify different curve characterizations. Since the arc-lenght parameter of $k-$principal direction curve of $\gamma$ is the same as the main curve $\gamma$, the impressions and expressions of theorems and definitions are obtained clearly and simply for $N_k-$slant helix. If the $N_k-$slant helix $\gamma_k$ has constant speed Darboux vector field, then $\gamma_k$ is called the curve of $N_k-$constant precession. It can be easily seen that the curve $\gamma$ is a slant helix for $k=0$, if $\gamma$ has constant speed Darboux then we obtain the characterization as the curve of constant precession.

**Acknowledgement.** The authors would like to express their sincere gratitude to the referees and Çetin Camcı for valuable suggestions to improve the paper.